\documentclass[12pt,reqno]{amsart}

\usepackage{setspace}
\usepackage{enumitem}
\usepackage{tikz-cd}
\usepackage{tikz}
\usetikzlibrary{positioning} 

\usepackage[linktocpage=true]{hyperref}

\makeatletter
\@namedef{subjclassname@2020}{%
  \textup{2020} Mathematics Subject Classification}
\makeatother
\newtheorem*{theorem*}{Theorem A}
\newtheorem*{theorem**}{Theorem B}

\usepackage[T1]{fontenc}
\usepackage{ae,aecompl}
\usepackage{graphics,cancel,graphicx}
\usepackage{epsfig,psfrag,manfnt}
\usepackage{mathrsfs}
\usepackage{graphicx,mathabx}
\usepackage{bbm,subfigure,rotating,bm,float,mathdots,wasysym,scalerel}

\newlength{\wdth}

\usepackage{stmaryrd}

\usepackage[all,cmtip]{xy}
\usepackage{amssymb,amsmath,amsfonts,amsthm,color}

\usepackage{scalerel,stackengine}
\stackMath
\newcommand\reallywidehat[1]{%
\savestack{\tmpbox}{\stretchto{%
  \scaleto{%
    \scalerel*[\widthof{\ensuremath{#1}}]{\kern-.6pt\bigwedge\kern-.6pt}%
    {\rule[-\textheight/2]{1ex}{\textheight}}
  }{\textheight}%
}{0.5ex}}%
\stackon[1pt]{#1}{\tmpbox}%
}

\newcommand{\mC}{\mathbb{C}}
\newcommand{\mD}{\mathbb{D}}

\newcommand{\mF}{\mathbb{F}}

\newcommand{\mN}{\mathbb{N}}

\newcommand{\mR}{\mathbb{R}}

\newcommand{\mZ}{\mathbb{Z}}

\newcommand{\blambda}{{\bm{\lambda}}}
\newcommand{\bmu}{{\bm{\mu}}}

\newtheorem{theorem}{Theorem}[section]
\newtheorem{lemma}[theorem]{Lemma}
\newtheorem{corollary}[theorem]{Corollary}

\theoremstyle{definition}
\newtheorem{remarks}[theorem]{Remarks}
\newtheorem{remark}[theorem]{Remark}

\theoremstyle{definition}

\theoremstyle{definition}

\theoremstyle{definition}

\begin{document}



\title[]{Maximal ideal space 
of some \\Banach algebras of Dirichlet series}

\author[]{Amol Sasane}
\address{Department of Mathematics \\London School of Economics\\
     Houghton Street\\ London WC2A 2AE\\ United Kingdom}
\email{A.J.Sasane@LSE.ac.uk}
 
\maketitle

\vspace{-1cm}
   
\begin{abstract}
\begin{spacing}{1.35}
Let $\mathscr{H}^\infty$ be the set of all Dirichlet series 
$f={\scaleobj{0.81}{\sum\limits_{n=1}^\infty}} \frac{a_n}{n^s}$ 
(where $a_n\in \mC$ for each $n$) that converge at each $s\in \mC_+$, 
such that $\|f\|_{\infty}:=\sup_{s\in \mC_+}|f(s)|<\infty$. 
Let $\mathscr{B}\subset \mathscr{H}^\infty$ 
be a Banach algebra containing the Dirichlet polynomials 
(Dirichlet series with finitely many nonzero terms) 
with a norm $\|\cdot\|_{{\scaleobj{0.81}{\mathscr{B}}}}$ 
such that the inclusion $\mathscr{B} \subset \mathscr{H}^\infty$ is continuous.
For $m\in \mN=\{1,2,3,\cdots\}$, let $\partial^{{\scaleobj{0.81}{-m}}}\mathscr{B}$ 
denote the Banach algebra consisting of all $f\in \mathscr{B}$ 
such that $f',\cdots,  f^{{\scaleobj{0.81}{(m)}}}\in \mathscr{B}$, 
with pointwise operations and the norm 
$
\textstyle 
\|f\|_{{\scaleobj{0.81}{\partial^{-m}\mathscr{B}}}}
\!=\!{\scaleobj{0.81}{\sum\limits_{\ell=0}^m}} 
\frac{1}{\ell!}\|f^{{\scaleobj{0.81}{(\ell)}}}\|_{{\scaleobj{0.81}{\mathscr{B}}}}.
$ 
Assuming that the Wiener $1/f$ property holds for $\mathscr{B}$ 
(that is, $\inf_{s\in \mC_+} |f(s)|>0$ implies $\frac{1}{f}\in \mathscr{B}$), 
it is shown that for all $m\in \mN$, 
the maximal ideal space $\text{M}(\partial^{{\scaleobj{0.81}{-m}}}\mathscr{B})$ 
of $\partial^{{\scaleobj{0.81}{-m}}}\mathscr{B}$ 
is homeomorphic to $\overline{\mD}^{{\scaleobj{0.81}{\mN}}}$, 
where $\overline{\mD}=\{z\in \mathbb{C}:|z|\le 1\}$. 
Examples of such Banach algebras are $\mathscr{H}^\infty$, 
the subalgebra $\mathscr{A}_u$ of $\mathscr{H}^\infty$ 
consisting of uniformly continuous functions in $\mC_+$, 
and the Wiener algebra $\mathscr{W}$ of Dirichlet series 
with $\|f\|_{{\scaleobj{0.81}{\mathscr{W}}}}
:={\scaleobj{0.81}{\sum\limits_{n=1}^\infty}} |a_n|<\infty$. 
Some consequences  (existence of logarithms, projective freeness, infinite Bass stable rank) 
are given as applications.
\end{spacing}
\end{abstract}

\vspace{-0.81cm}

\section{Introduction}

\begin{spacing}{1.1}
\noindent The aim of this article is 
to determine the maximal ideal space 
of a particular family $\{\partial^{-m} \mathscr{B}\}_{m\in \mN}$ (defined  below) of Banach algebras 
that are contained in the Hardy algebra $\mathscr{H}^\infty$ of Dirichlet series. 
\end{spacing}

\begin{spacing}{1.3}
Set $\mC_+:=\{s\in \mC: \text{Re}(s)>0\}$. 
Let $\mathscr{H}^\infty$ be the set of all Dirichlet series 
$f={\scaleobj{0.72}{\sum\limits_{n=1}^\infty\cfrac{a_n}{n^s}}}$, 
where $a_n\in \mC$ for each $n\in \mN$, that converge for all $s\in \mC_+$, 
such that $\|f\|_{\infty}:=\sup_{s\in \mC_+}|f(s)|<\infty$. 
With pointwise operations and the supremum norm, 
$\mathscr{H}^\infty$ is a Banach algebra.
\end{spacing}

\subsection*{The Banach algebras $\partial^{-m}\mathscr{B}$} 

Throughout, $\mathscr{B}\subset \mathscr{H}^\infty$ will denote a Banach algebra 
with a norm $\|\cdot\|_{{\scaleobj{0.81}{\mathscr{B}}}}$. 
For $m\in \mN=\{1,2,3,\cdots\}$, let $\partial^{{\scaleobj{0.81}{-m}}}\mathscr{B}$ 
denote the Banach algebra consisting of all $f\in \mathscr{B}$ such that $f',\cdots,  f^{{\scaleobj{0.81}{(m)}}}\in \mathscr{B}$, with pointwise operations and the norm 
 $$
\textstyle 
\|f\|_{{\scaleobj{0.81}{\partial^{-m}\mathscr{B}}}}\!=\!{\scaleobj{0.81}{\sum\limits_{\ell=0}^m}} \;\!
 {\scaleobj{0.81}{\cfrac{1}{\ell!}}} \;\!\|f^{{\scaleobj{0.81}{(\ell)}}}\|_{{\scaleobj{0.81}{\mathscr{B}}}}.
$$ 
\begin{spacing}{1.35}
\noindent We note that if $f={\scaleobj{0.81}{\sum\limits_{n=1}^\infty}} {\scaleobj{0.72}{\cfrac{a_n}{n^s}}}$ converges for all $s\in \mC_+$, then it converges uniformly on compact sets contained in $\mC_+$, and hence by Weierstrass's theorem on uniform limits of holomorphic functions, $f^{(\ell)}$ is obtained by termwise differentiation, so that for all  $\ell\in \mN$, we have 
\end{spacing}
\vspace{-0.45cm}
$$
\textstyle 
f^{(\ell)}={\scaleobj{0.81}{\sum\limits_{n=1}^\infty}}(-1)^\ell (\log n)^\ell {\scaleobj{0.81}{\cfrac{a_n}{n^s}}}\text{ in }\mC_+. 
$$
Let $\mathscr{P}$ denote the set of {\em Dirichlet polynomials}, 
that is, Dirichlet series with finite support, 
$$
\textstyle 
\mathscr{P} = \big\{ p={\scaleobj{0.81}{\sum\limits_{n=1}^N}} {\scaleobj{0.72}{\cfrac{a_n}{n^s}}}:  
N\in \mN, \;a_1,\cdots, a_N\in \mC \big\} \subset \mathscr{H}^\infty.
$$
Let $\mathscr{A}_u$ be the subset of $\mathscr{H}^\infty$ 
of Dirichlet series that are uniformly continuous in $\mC_+$.
Another description of $\mathscr{A}_u$ is that 
it is the closure of Dirichlet polynomials in the $\|\cdot\|_\infty$-norm, 
see, e.g., \cite[Theorem~2.3]{ABG}.

Let $\mathscr{W}$ denote the set of all Dirichlet series $f={\scaleobj{0.81}{\sum\limits_{n=1}^\infty}}  {\scaleobj{0.72}{\cfrac{a_n}{n^s}}}$ such that 
$$
\textstyle \|f\|_{1}:={\scaleobj{0.81}{\sum\limits_{n=1}^\infty}} |a_n|<\infty.
$$ 
With pointwise operations and the $\|\cdot\|_1$ norm, $\mathscr{W}$ is a Banach algebra. It is clear that 
$$
\mathscr{W} \subset \mathscr{A}_u \subset \mathscr{H}^\infty.
$$
In the case of $\mathscr{W}$, an analogue of the classical Wiener $1/f$ lemma (\cite[p.91]{Wie}) for the unit circle holds, 
that is, if $f\in \mathscr{W}$ is such that $\inf_{s\in \mC_+} |f(s)|>0$, then $\frac{1}{f}\in \mathscr{W}$ 
(see, e.g., \cite[Thm.~1]{HewWil}, and also \cite{GooNew} for an elementary proof). 

We say that a Banach algebra $\mathscr{B}\subset \mathscr{H}^\infty$ has the {\em Wiener property}  if 
$$
\textstyle 
\text{(}W\text{)} 
\;\;\text{For all } f\in \mathscr{B} \text{ satisfying } \inf_{s\in \mC_+} |f(s)|>0, 
\text{ we have } \frac{1}{f} \in \mathscr{B}.
$$
The Banach algebra $\mathscr{H}^\infty$ also possesses the Wiener property $(W)$ (see, e.g., \cite[Theorem~2.6]{Bon}).

\begin{lemma}
 $\mathscr{A}_u$ possesses the Wiener property $(W)$. 
\end{lemma}
\begin{proof}  Let $f\in \mathscr{A}_u$ satisfy $d:=\inf_{s\in \mC_+} |f(s)|>0$. 
As $\mathscr{A}_u\subset \mathscr{H}^\infty$, it follows that $\frac{1}{f}\in \mathscr{H}^\infty$. 
Moreover, $\frac{1}{f}$ is uniformly continuous in $\mC_+$: for all $z,w\in \mC_+$, 
we have 
$$
\textstyle |\frac{1}{f}(w)-\frac{1}{f}(z)|=\frac{|f(z)-f(w)|}{|f(z)||f(w)|}\le \frac{1}{d^2} |f(w)-f(z)|,
$$
 and $f$ is uniformly continuous in $\mC_+$.
\end{proof}

\noindent Let $A$ be a commutative unital complex semisimple Banach algebra. The dual space $A^*$ of $A$ consists of all continuous linear complex-valued maps defined on $A$. The {\em maximal ideal
space} $\text{M}(A)$ of $A$  is the set of all nonzero multiplicative elements in $A^*$ (the kernels of which are then in one-to-one correspondence with the maximal ideals of $A$). As $\text{M}(A)$ is a subset of $A^*$, it inherits the  weak-$\ast$ topology of $A^*$, called the {\em Gelfand topology} on $\text{M}(A)$. 
 The topological space $\text{M}(A)$ is a compact Hausdorff space, and is contained in the unit sphere of  the Banach space $A^*$  with the operator norm, $\|\varphi\|=\sup_{a\in A, \;\|a\|\le 1} |\varphi(a)|$ for all $\varphi \in A^*$.  
 Let $C(\text{M}(A))$ denote the Banach algebra  of 
complex-valued continuous functions on $\text{M}(A)$ with pointwise operations and the supremum norm. 
 The {\em Gelfand transform} $\widehat{a}\in C(\text{M}(A))$ of an element $a\in A$ is  
 defined by $\widehat{a}(\varphi):=\varphi(a)$ for all $\varphi \in \text{M}(A)$. 
 
\subsection*{Main result} 

The main result in this article is the following.

\begin{theorem} 
\label{main_theorem}
Let $m\in \mN,$ and let the Banach algebra $\mathscr{B}$ be such that 
\begin{itemize}
\item[${\scaleobj{0.81}{\bullet}}$] $\mathscr{P}\subset \mathscr{B}\subset \mathscr{H}^\infty$
\item[${\scaleobj{0.81}{\bullet}}$] there exists a $C>0$ such that for all $f\in \mathscr{B},$ 
$\|f\|_\infty \le C\|f\|_{{\scaleobj{0.81}{\mathscr{B}}}}$
 \item[${\scaleobj{0.81}{\bullet}}$] $\mathscr{B}$ possess the Wiener property $(W)$. 
\end{itemize}
Then the maximal ideal space of $\partial^{-m} \mathscr{B}$ is homeomorphic to $\overline{\mD}^{{\scaleobj{0.81}{\mN}}}$.
\end{theorem}

\noindent Here each factor $\overline{\mD}$ has the usual Euclidean topology inherited from $\mC$, and $\overline{\mD}^{{\scaleobj{0.81}{\mN}}}$ is given the product topology.

\subsection*{Examples}
 Examples of such Banach algebras $\mathscr{B}$ are $\mathscr{H}^\infty, \mathscr{A}_u$ and $\mathscr{W}$. 
Given a subset $S\subset i\mR$, the Banach algebra 
$$
\textstyle \mathscr{H}^\infty_S:=\{ f\in \mathscr{H}^\infty: f \text{ has a continuous extension to }S\},
$$
with pointwise operations and the norm $\|\cdot\|_\infty$, is also one that satisfies the assumptions of Theorem~\ref{main_theorem}. The Wiener property $(W)$ for $\mathscr{H}^\infty_S$ is an immediate consequence of that for $\mathscr{H}^\infty$.

\subsection*{Organisation of the article} 

In Section~\ref{section_2}, we will prove Theorem~\ref{main_theorem}, and in Section~\ref{section_3}, some corollaries  (existence of logarithms, projective freeness, infinite Bass stable rank) are given as applications.

\section{Proof of  the main result}
\label{section_2}

\noindent We first show the following, which will be used to prove Theorem~\ref{main_theorem}. 

\begin{lemma}
\label{lemma_10_3_25_1548}
If $m\in \mN,$ then $\partial^{-m}\mathscr{B}\subset \mathscr{A}_u$. 
\end{lemma}
\begin{proof} Let $f\in \partial^{-m}\mathscr{B}$. As $m\geq 1$, $f'\in \mathscr{B}\subset \mathscr{H}^\infty$. For $z,w\in \mC_+$, let $[z,w]$ denote the straight line segment joining $z$ to $w$. By the fundamental theorem of contour integration, $f(w)-f(z)=\int_{[z,w]} f'(\zeta)d\zeta$. By the $ML$-inequality, 
$$
\textstyle 
|f(w)-f(z)|\le |w-z| \max\limits_{\zeta\in [z,w]} |f'(\zeta)|\le |w-z|\|f'\|_\infty.
$$
  Thus $f$ is uniformly continuous in $\mC_+$. Also $f\in \mathscr{H}^\infty$. So $f\in  \mathscr{A}_u$. 
\end{proof}

\noindent Let $p_1<p_2<p_3<\cdots$ be the sequence of all primes arranged in increasing order. By the fundamental theorem of arithmetic, every $n\in \mN$ 
can be written uniquely in the form 
$$
\textstyle 
n={\scaleobj{0.81}{\prod\limits_{k=1}^\infty}} \;p_k^{{\scaleobj{0.81}{\nu_{p_k}\!(n)}}},
$$ 
where $\nu_{p_k}(n)\in \mN\cup\{0\}$ denotes the largest integer 
$m$ such that $p_k^m$ divides $n$. 

\begin{proof}[Proof of Theorem~\ref{main_theorem}.]
For $\blambda \!=\!(\lambda_1,\lambda_2,\cdots)\!\in \overline{\mD}^{{\scaleobj{0.81}{\mN}}}$, define  $\varphi_{\blambda}\!:\! \mathscr{P}\!\to \mC$ by 
$$
\textstyle 
\varphi_{\blambda}(p)={\scaleobj{0.81}{\sum\limits_{n=1}^N}}\;\! a_n \;\!{\scaleobj{0.81}{\prod\limits_{k=1}^\infty}} \;\!\lambda_k^{{\scaleobj{0.81}{\nu_{p_k}\!(n)}}},
\;\text{ for }p={\scaleobj{0.81}{\sum\limits_{n=1}^N}} \;\!{\scaleobj{0.81}{\cfrac{a_n}{n^s}}}\in \mathscr{P}.
$$
For each $n\in \mN_*$, since 
$$
\textstyle 
|{\scaleobj{0.81}{\prod\limits_{k=1}^\infty}}\; \lambda_k^{{\scaleobj{0.81}{\nu_{p_k}\!(n)}}}|\le 1
$$ 
we have that $|\varphi_{\blambda}(p)|\le \|p\|_1\le \|p\|_\infty$.  
 As $m\ge 1$, it follows from Lemma~\ref{lemma_10_3_25_1548} that $\mathscr{B}\subset \mathscr{A}_u$. 
 So $\mathscr{P}$ is dense in $\partial^{-m}\mathscr{B}$ in the $\|\cdot\|_\infty$-norm.  Given $f\in \partial^{-m}\mathscr{B}$, 
 let $(p_n)_{n\in \mN}$ be a sequence in $\mathscr{P}$ that converges to $f$ in the  $\|\cdot\|_\infty$-norm. 
 Then $(\varphi_\blambda (p_n))_{n\in \mN}$ is a Cauchy sequence in $\mC$ (as $|\varphi_\blambda (p_n)-\varphi_\blambda(p_m)|\le \|p_n-p_m\|_\infty$), and hence convergent. Define 
 $$
 \textstyle \varphi_\blambda(f)=\lim\limits_{n\to \infty} \varphi_\blambda (p_n).
 $$
 Then $\varphi_\blambda: \partial^{-m} \mathscr{B}\to \mC$ is well-defined: if $(\widetilde{p}_n)_{n\in \mN}$ is another sequence of approximating Dirichlet polynomials, then 
 $$
 \textstyle |\varphi_\blambda (\widetilde{p}_n)-\varphi_\blambda(p_n)|\le \| \widetilde{p}_n-p_n\|_\infty 
 \le \|\widetilde{p}_n-f\|_\infty +\|f-p_n\|_\infty \to 0 \text{ as } n\to \infty,
 $$
 and so $\lim\limits_{n\to \infty} \varphi_\blambda (\widetilde{p}_n)=\lim\limits_{n\to \infty} \varphi_\blambda (p_n)+\lim\limits_{n\to \infty} \varphi_\blambda (\widetilde{p}_n-p_n)=\lim\limits_{n\to \infty} \varphi_\blambda (p_n)+0$.
 
 \goodbreak 
 
 \noindent 
We claim that the map $\varphi_\blambda$ is a complex homomorphism. 
It is enough to show linearity and multiplicativity on $\mathscr{P}$, since it then extends to $\partial^{-m} \mathscr{B}$ by the algebra of limits and the continuity of addition, scalar multiplication and multiplication on $\mathscr{P}$ in the $\|\cdot\|_\infty$-norm.  Linearity is clear, so we just show  multiplicativity: 
$$
\textstyle 
\begin{array}{rcl}
\varphi_\blambda (pq)
\!\!\!&=&\!\!\!
{\scaleobj{0.81}{\sum\limits_{n=1}^N}} \;({\scaleobj{0.81}{\sum\limits_{d\vert n} }} \; a_d b_{\frac{n}{d}} )
{\scaleobj{0.81}{\prod\limits_{k=1}^\infty}} \;\lambda_k^{{\scaleobj{0.81}{\nu_{p_k}\!(n)}}}
=
{\scaleobj{0.81}{\sum\limits_{n=1}^N}}\; ({\scaleobj{0.81}{\sum\limits_{d\vert n} }} \;a_d b_{\frac{n}{d}} )
{\scaleobj{0.81}{\prod\limits_{k=1}^\infty}} \;\lambda_k^{{\scaleobj{0.81}{\nu_{p_k}\!(d) + \nu_{p_k}\!(\frac{n}{d})}}} \\
\!\!\!&=&\!\!\! 
({\scaleobj{0.81}{\sum\limits_{d=1}^N}} \;a_d {\scaleobj{0.81}{\prod\limits_{k=1}^\infty}} \;\lambda_k^{{\scaleobj{0.81}{\nu_{p_k}\!(d) }}} )
({\scaleobj{0.81}{\sum\limits_{\widetilde{d}=1}^N}} \;b_{\widetilde{d}} {\scaleobj{0.81}{\prod\limits_{k=1}^\infty}}\; \lambda_k^{{\scaleobj{0.81}{\nu_{p_k}\!(\widetilde{d}) }}} )
=
\varphi_\blambda (p)\;\!\varphi_\blambda (q)
\end{array}
$$
for all $p={\scaleobj{0.81}{\sum\limits_{n=0}^N}}\;\! {\scaleobj{0.72}{\cfrac{a_n}{n^s}}}, \;q={\scaleobj{0.81}{\sum\limits_{n=0}^N}}\;\!{\scaleobj{0.72}{ \cfrac{b_n}{n^s}}}\in 
\mathscr{P}$. Finally, $\varphi_\blambda$ is bounded, because
$$
\textstyle 
\begin{array}{rcl}
|\varphi_\blambda(f)|=\big|\lim\limits_{n\to \infty} \varphi_\blambda (p_n)|=\lim\limits_{n\to \infty}|\varphi_\blambda (p_n)|
\!\!\!&\le&\!\!\! \lim\limits_{n\to \infty}\| p_n\|_\infty =\|f\|_\infty\\[0.27cm]
\!\!\!&\le&\!\!\!  C\|f\|_{{\scaleobj{0.81}{\mathscr{B}}}} \le  C\|f\|_{{\scaleobj{0.81}{\partial^{-m}\mathscr{B}}}},
\end{array}
$$
where $f\in \partial^{-m}\mathscr{B}$, and $(p_n)_{n\in \mN}$ is an approximating sequence in $\mathscr{P}$ for $f$ in the $\|\cdot\|_\infty$-norm. Note that in particular, we have $\|\varphi_\blambda(f)|\le \|f\|_\infty$. 

Let $\blambda=(\lambda_1, \lambda_2,\cdots) ,\bmu=(\mu_1,\mu_2,\cdots)$ be distinct elements of $\overline{\mD}^{{\scaleobj{0.81}{\mN}}}$. Then there exists an $k_{{\scaleobj{0.81}{ *}}}\in \mN_*$ such that $\lambda_{k_{{\scaleobj{0.81}{ *}}}}\neq \mu_{k_{{\scaleobj{0.81}{ *}}}}$. 
As $\partial^{-m} \mathscr{B}$ contains $\frac{1}{2^s}, \frac{1}{3^s},\cdots\in \mathscr{P}$, we have 
$$
\textstyle 
\varphi_\blambda({\scaleobj{0.72}{\cfrac{1}{p_{k_{{\scaleobj{0.81}{ *}}}}^s}}})=1\cdot \lambda_1^0 \cdots\lambda_{k_{{\scaleobj{0.81}{ *}}}-1}^0 \lambda_{k_{{\scaleobj{0.81}{ *}}}}^1 \lambda_{k_{{\scaleobj{0.81}{ *}}}+1}^0 \cdots =\lambda_{k_{{\scaleobj{0.81}{ *}}}}\neq \mu_{k_{{\scaleobj{0.81}{ *}}}}
 =\varphi_\bmu({\scaleobj{0.72}{\cfrac{1}{p_{k_{{\scaleobj{0.81}{ *}}}}^s}}}).
$$
Thus $\blambda \mapsto \varphi_\blambda$ embeds $\overline{\mD}^{{\scaleobj{0.81}{\mN}}}$ in the maximal ideal space of $\partial^{-m} \mathscr{B}$. 

\begin{spacing}{1.15}
We claim  that the inclusion $\overline{\mD}^{{\scaleobj{0.81}{\mN}}}\subset \text{M}(\partial^{-m} \mathscr{B})$ is continuous. 
Let $(\blambda_i )_{i\in I}$ be a net in $\overline{\mD}^{{\scaleobj{0.81}{\mN}}}$ which is convergent to $\blambda \in \overline{\mD}^{{\scaleobj{0.81}{\mN}}}$. Let $\epsilon>0$ and $ f\in 
\partial^{-m}\mathscr{B}.
$ 
 Then there exists a 
 \end{spacing}
 \vspace{-0.45cm}
 $$
 \textstyle p= {\scaleobj{0.81}{\sum\limits_{n=0}^N}}\; {\scaleobj{0.72}{\cfrac{b_n}{n^s}}}\in \mathscr{P},
 $$
 \noindent 
   such that $\|f-p\|_\infty<\frac{\epsilon}{4}$. 
 Let $p_1,\cdots, p_{k_N}$ be the only primes which appear in the prime factorisation of $1,\cdots, N$. 
 If $\succcurlyeq$ denotes the order on the directed set $I$, then there exists an $i_*\in I$ such that for all $i\succcurlyeq i_*$, 
 $$
 \textstyle 
 {\scaleobj{0.81}{\sum\limits_{n=1}^N}} \; |b_n| \;\! \big|{\scaleobj{0.81}{\prod\limits_{k=1}^{k_N}}}\; \lambda_{i,k}^{{\scaleobj{0.81}{\nu_{p_k}\!(n) }}}
 \!-\!{\scaleobj{0.81}{\prod\limits_{k=1}^{k_N}}}\; \lambda_{i_*,k}^{{\scaleobj{0.81}{\nu_{p_k}\!(n) }}}\big|<\frac{\epsilon}{2}, 
 $$
 and so 
 $|\varphi_{\blambda_i}(p)-\varphi_{\blambda}(p)| \le  {\scaleobj{0.81}{\sum\limits_{n=1}^N}} \; |b_n| \;\! \big|{\scaleobj{0.81}{\prod\limits_{k=1}^{k_N}}}\; \lambda_{i,k}^{{\scaleobj{0.81}{\nu_{p_k}\!(n) }}}
 \!-\!{\scaleobj{0.81}{\prod\limits_{k=1}^{k_N}}}\; \lambda_{i_*,k}^{{\scaleobj{0.81}{\nu_{p_k}\!(n) }}}\big| < \frac{\epsilon}{2}$. 
 Thus 
 $$
 \textstyle 
 \begin{array}{rcl}
 |\varphi_{\blambda_i}(f)-\varphi_{\blambda}(f)|
 \!\!\!&\le&\!\!\! 
  |\varphi_{\blambda_i}(p)-\varphi_{\blambda}(p)|+|\varphi_{\blambda_i}(f-p)|+|\varphi_{\blambda}(f-p)|
  \\[0.2cm]
  \!\!\!&\le&\!\!\!  \frac{\epsilon}{2}+ \|f-p\|_\infty + \|f-p\|_\infty \le \frac{\epsilon}{2}+\frac{\epsilon}{4}+\frac{\epsilon}{4}=\epsilon
  \end{array}
  $$
 for all $i\succcurlyeq i_*$. 
 Hence $(\varphi_{\blambda_i})_{i\in I}$ converges to $\varphi_{\blambda}$ in the weak-$\ast$ topology, i.e., the Gelfand topology on  the maximal ideal space of $\partial^{-m}\mathscr{B}$.

\noindent  Next we will show that every complex homomorphism is of the form $\varphi_\blambda$ for some $\blambda \in \overline{\mD}^{{\scaleobj{0.81}{\mN}}}$. 
 
  Let $\varphi\in M(\partial^{-m}\mathscr{B})$. Define 
  $$
 \textstyle 
 \blambda=\pmb{(}\;\!\varphi({\scaleobj{0.72}{\cfrac{1}{2^s}}}), \varphi({\scaleobj{0.72}{\cfrac{1}{3^s}}}), \varphi({\scaleobj{0.72}{\cfrac{1}{5^s}}}), \cdots\pmb{)}.
 $$ 
 We first show that for all $f\in \partial^{-m}\mathscr{B}$, we have 
 $$
 \textstyle 
\quad \quad  \quad \quad \quad\quad \quad \quad \quad\quad \;\;
 |\varphi(f)|\le \|f\|_\infty.
\quad \quad \quad \quad \quad \quad \quad \quad\quad (\asterisk)
 $$
 Suppose first that $f$ also satisfies 
 $$
 \textstyle 
 \inf\limits_{s\in \mC_+} |f(s)|>0.
 $$
 As $\mathscr{B}$ possesses the Wiener property $(W)$, we have $\frac{1}{f}\in \mathscr{B}$. 
 Differentiating, we get successively that 
 $$
 \textstyle 
 ({\scaleobj{0.72}{\cfrac{1}{f}}})'=-{\scaleobj{0.72}{\cfrac{f'}{f^2}}},\quad ({\scaleobj{0.72}{\cfrac{1}{f}}})''=-{\scaleobj{0.72}{\cfrac{f''f^2-2f (f')^2}{f^4}}},\quad \cdots,
 $$
 \begin{spacing}{1.15}
 \noindent and so (since $f,f',\cdots, f^{(m)}\in \mathscr{B}$), we conclude that $\frac{1}{f}\in \partial^{-m}\mathscr{B}$. So we have shown that if $0$ does not belong to the closure of the range of $f\in \partial^{-m}\mathscr{B}$, then $\frac{1}{f}\in  \partial^{-m}\mathscr{B}$, 
 and in particular $1=\varphi(1)=\varphi(f)\varphi(\frac{1}{f})$, showing that $\varphi(f)\neq 0$. 
 Replacing $f$ by $f-c$, where $c\in \mC$, we conclude that if $c$ does not belong to the closure of the range of $f$, then $\varphi(f)\neq c$. Thus $\varphi(f)$ belongs to the closure of the range of $f$. In particular, $|\varphi(f)|\leq \|f\|_\infty$, as wanted. 
 \end{spacing}

 Applying this to $f:=\frac{1}{p_k^s}$ yields $|\lambda_k|\le 1$, $k\in \mN_*$, and so $\blambda \in \overline{\mD}^{{\scaleobj{0.81}{\mN}}}$. 
 
 Since $\partial^{-m}\mathscr{B} \subset \mathscr{A}_u$, any 
 $f \in \partial^{-m}\mathscr{B}
 $ 
 can be approximated in the $\|\cdot\|_\infty$-norm by a sequence $(p_n)_{n\in \mN}$ of Dirichlet polynomials. 
 But ($\asterisk$) shows that $\varphi$ is continuous in the $\|\cdot\|_\infty$-norm, giving 
 $$
 \textstyle 
 \varphi(f)
 \!=\! \varphi(\lim\limits_{n\to \infty} p_n)
 \!=\! \lim\limits_{n\to \infty}\varphi(p_n)
 \!=\! \lim\limits_{n\to \infty}\varphi_\blambda (p_n)
 \!=\! \varphi_\blambda (\lim\limits_{n\to \infty} p_n)
 \!=\! \varphi_\blambda (f).
 $$

 \begin{spacing}{1.15}
\noindent We have seen that the Gelfand topology of the maximal ideal space of $\partial^{-m}\mathscr{B}$ is weaker/coarser than the product topology of $\overline{\mD}^{{\scaleobj{0.81}{\mN}}}$.  
As the Gelfand topology is Hausdorff, and $\overline{\mD}^{{\scaleobj{0.81}{\mN}}}$ is compact (Tychonoff's theorem), the two topologies coincide (see, e.g., \cite[14, \S3.8]{Rud}, stating that if $\tau_1\subset \tau_2$ are topologies on a set $X$, such that $\tau_1$ is Hausdorff and $\tau_2$ is compact, then $\tau_1=\tau_2$). 
\end{spacing}
\vspace{-0.45cm}
\end{proof}

\begin{remark} 
\begin{spacing}{1.1}
The theorem and its proof above also works for $m=0$ if $\mathscr{B}$ is $\mathscr{A}_u$ or $\mathscr{W}$. 
The maximal ideal space of $\mathscr{W}$ as being homeomorphic to $\overline{\mD}^{{\scaleobj{0.81}{\mN}}}$ was shown in \cite{Wla}. 
\end{spacing}
\end{remark}

\vspace{-0.6cm}

\section{Some consequences}
\label{section_3}

\noindent Throughout this section, we will assume that 
 $m\in \mN,$ and  $\mathscr{B}$ is a Banach algebra such that 
\begin{itemize}
\item[${\scaleobj{0.81}{\bullet}}$] $\mathscr{P}\subset \mathscr{B}\subset \mathscr{H}^\infty$
\item[${\scaleobj{0.81}{\bullet}}$] there exists a $C>0$ such that for all $f\in \mathscr{B},$ 
$\|f\|_\infty \le C\|f\|_{{\scaleobj{0.81}{\mathscr{B}}}}$
 \item[${\scaleobj{0.81}{\bullet}}$] $\mathscr{B}$ possess the Wiener property $(W)$. 
\end{itemize}

\subsection*{Contractibility of $\text{M}(\partial^{-m}\mathscr{B})$} 
Recall that a topological space $X$ is {\em contractible} 
if the identity map $\textrm{id}_X:X \to X$ is null-homotopic, 
i.e., there exist an element $x_*\in X$ and a continuous map 
$\textrm{H}:[0,1]\times X \to X$ such that $\textrm{H}(0, \cdot)=\textrm{id}_X$ 
and $\textrm{H}(1,x)= x_*$ for all $x\in X$.

\begin{corollary}
$\text{\em M}(\partial^{-m}\mathscr{B})$ is contractible.
\end{corollary}
\begin{proof} \begin{spacing}{1.1}It suffices to show $\overline{\mD}^{{\scaleobj{0.81}{\mN}}}$ is contractible. Let ${\bm{x}}_*\!=\!\mathbf{0}\!=\!\pmb{(}0,0,\cdots\pmb{)}
\!\in\! \overline{\mD}^{{\scaleobj{0.81}{\mN}}}$, and $\textrm{H}(t,{\bm{x}})\!=\!(1-t){\bm{x}}\!=\!\pmb{(}(1-t)x_1, (1-t)x_2,\cdots\pmb{)}$ for ${\bm{x}}=\pmb{(}x_1,x_2,\cdots\pmb{)}\in   \overline{\mD}^{{\scaleobj{0.81}{\mN}}}$ and $t\in [0,1]$. Then $\textrm{H}$ is continuous, $\textrm{H}(0, \cdot)=\textrm{id}_X$, and $\textrm{H}(1,{\bm{x}})= {\bm{x}}_*$ for all ${\bm{x}}\in \overline{\mD}^{{\scaleobj{0.81}{\mN}}}$. 
\end{spacing}
\vspace{-0.45cm}
\end{proof}

\subsection*{Existence of logarithms}

For a unital commutative complex Banach algebra $A$, 
 the multiplicative group of all invertible elements of $A$ is denoted by $A^{-1}$. 
 Then  $e^A:=\{e^a:a\in A\}$  is a subgroup of $A^{-1}$.
By the Arens-Royden theorem (see, e.g., \cite[Theorem, p.295]{Roy}), 
the group $A^{-1}/e^{A}$ is isomorphic 
to the first \v{C}ech cohomology group $H^1(\text{M}(A), \mZ)$ of $\text{M}(A)$ with integer coefficients.  
For background on \v{C}ech cohomology, see, e.g., \cite{EilSte}. 
For a contractible space, 
all cohomology groups are trivial (see, e.g., \cite[IX,\,Theorem~3.4]{EilSte}). 

\begin{corollary}
$(\partial^{-m}\mathscr{B})^{-1}=e^{\partial^{-m}\mathscr{B}}$.
\end{corollary}

\subsection*{Projective freeness} 
For a commutative unital ring $A$ with unit element denoted by $1$, $A^{n\times n}$ denotes 
the $n \times n$ matrix ring over $A$, and 
$\text{GL}_n(A) \subset A^{n\times n}$ denotes the group of invertible matrices. 
A commutative unital ring $A$ is {\em projective free} if 
every finitely generated projective $A$-module is free. 
If $A$-modules $M,N$ are isomorphic, then we write $M\cong N$. 
If  $M$ is a finitely generated $A$-module, then  
(i) $M$ is {\em free} if $M \cong A^k$ for some  $k\in \mN\cup\{ 0\}$, and 
(ii) $M$ is {\em projective} if there exists an $A$-module $N$ 
and an  $n\in \mN\cup\{ 0\}$ such that $M\oplus N \cong A^n$. 
In terms of matrices (see, e.g., \cite[Proposition~2.6]{Coh}), 
the ring $A$ is projective free if and only if 
every idempotent matrix $P$ is conjugate (by an invertible matrix $S$) 
to a diagonal matrix with elements $1$ and $0$ on the diagonal, 
 i.e., for all $n\in \mN$ and every $P\in A^{n\times n}$ satisfying $P^2=P$, 
there exists an $S\in \text{GL}_n(A)$ such that  for some  $k\in \mN\cup\{ 0\}$, 
 $
S^{-1} P S
=
[\begin{smallmatrix}
I_k & 0\\
0 & 0
\end{smallmatrix}].
$ 

\noindent 
In 1976, it was shown independently by Quillen and Suslin that 
if $\mF$ is a field, then the polynomial ring $\mF[x_1, \dots , x_n]$ 
is projective free, settling Serre's conjecture from 1955 (see \cite{Lam}). 
In the context of a commutative  unital complex Banach algebra $A$, 
\cite[Theorem~4.1]{BruSas23} 
says that the contractibility of the maximal ideal space $\text{M}(A)$ 
is sufficient for $A$ to be projective free.

\begin{corollary}
$\partial^{-m}\mathscr{B}$ is a projective free ring. 
\end{corollary}

\subsection*{Bass stable rank}

In algebraic $K$-theory, the notion of stable rank of a ring
was introduced to facilitate $K$-theoretic computations 
\cite{Bas}.
Let $A$ be a unital commutative ring with unit element denoted by $1$. 
An element $(a_1,\cdots, a_n)\!\in\! A^n$ is {\em unimodular} if
there exist  $b_1,\cdots, b_n\in A$ such that
  $
 b_1 a_1+\cdots +b_n a_n=1. 
 $
The set of all unimodular elements of $A^n$ is denoted by $U_n (A)$. 
 We call $(a_1,\cdots, a_{n+1})\in U_{n+1}(A)
$ 
 {\em reducible} if there exist  $x_1,\cdots, x_n\!\in \!A $ such
that
 $
(a_1\!+\!x_1 a_{n+1},\;\!\cdots, \;\! a_n\!+\! x_n a_{n+1})\!\in\! U_n(A).
$ 
The {\em Bass stable rank} of $A$ is the least  $n\in \mN$ for
which every element in $ U_{n+1}(A)$ is reducible. The {\em Bass stable rank of $A$ is infinite} if there is no such
 $n$.
Analogous to the result that the Bass stable
rank of the  infinite polydisc algebra $A(\mD^\infty)$ is infinite (see, e.g.,  \cite[Proposition~1]{Mor}), we show the following.

\begin{corollary}
\label{Bass_stable_rank_of_A_DS}
 The Bass stable rank of $\partial^{-m}\mathscr{B}$ is infinite.
\end{corollary}
\begin{proof} Fix $n\in \mN$. Let $f_1,\cdots, f_{n+1}\in \mathscr{P} \subset \partial^{-m}\mathscr{B}$ be given by 
$$
\textstyle 
f_1\!=\!{\scaleobj{0.72}{\cfrac{1}{2^s}}},\;\;\cdots,\;\; f_n\!=\! {\scaleobj{0.72}{\cfrac{1}{p_n^s}}}, \;\;f_{n+1}\!=\! {\scaleobj{0.81}{\prod\limits_{j=1}^n}}\;\! \big(1\!-\!{\scaleobj{0.72}{\cfrac{1}{(p_j p_{n+j})^s}}} \big) .
$$
Then $(f_1,\cdots, f_{n+1})\in U_{n+1}( \partial^{-m}\mathscr{B})$ because 
by expanding the product defining $f_{n+1}$, we
obtain
$$
\textstyle 
f_{n+1}=1-{\scaleobj{0.72}{\cfrac{1}{2^s}}}\cdot  g_1-\cdots -{\scaleobj{0.72}{\cfrac{1}{p_n^s}}}\cdot g_n=1-f_1 g_1-\cdots-f_n g_n,
$$
for suitably defined $g_1,\cdots, g_n\in  \mathscr{P} \subset \partial^{-m}\mathscr{B}$, and so with $g_{n+1}:=1$, we get 
$f_1g_1+\cdots+f_n g_n+f_{n+1} g_{n+1}=1$. 
 Let
$(f_1,\cdots, f_{n+1})$ be reducible, and  $x_1,\cdots,x_n\in
\partial^{-m}\mathscr{B}$ be such that
$$
\textstyle \big(\;\!{\scaleobj{0.72}{\cfrac{1}{2^s}}}+x_1f_{n+1},\;\cdots,\;{\scaleobj{0.72}{\cfrac{1}{p_n^s}}}+x_nf_{n+1}\;\!\big)\in U_n(\partial^{-m}\mathscr{B}).
$$ 
Let $y_1,\cdots , y_n\in \partial^{-m}\mathscr{B}$ be such that 
$$
\textstyle 
\big(\;\! {\scaleobj{0.72}{\cfrac{1}{2^s}}}+x_1f_{n+1}\big)y_1+\cdots
+\big(\;\!{\scaleobj{0.72}{\cfrac{1}{p_n^s}}}+x_nf_{n+1}\big)y_n=1.
$$
Taking the Gelfand transform,  and denoting the variable in the infinite polydisc $\overline{\mD}^{{\scaleobj{0.81}{\mN}}}$ by  
$\bm{z}=(z_1,z_2,z_3,\cdots)$, we obtain
$$
\textstyle
\quad\quad\quad\quad\quad\;\;
(z_1+ \widehat{x}_1\widehat{f}_{n+1} )\widehat{y}_1+\cdots+ (z_n+\widehat{x}_n\widehat{f}_{n+1})\widehat{y}_n=1 .
\quad\quad\quad\quad\quad\;\; (\star)
$$
Let ${\bm{x}}:=(\widehat{x}_1,\cdots, \widehat{x}_n)$.  For
${\bm{z}}=(z_1,\cdots, z_n)\in \mC^n$, we define
$$
\textstyle 
\mathbf{\Phi}({\bm{z}})\!=\!\bigg\{\!
{\scaleobj{0.9}{\begin{array}{l} 
 -{\bm{x}}(z_1,\cdots, z_n,\overline{z_1},\cdots, \overline{z_n},0,\cdots) 
{\scaleobj{0.81}{\prod\limits_{j=1}^n}} \;\!(1\!-\!|z_j|^2)  
\textrm{ if } |z_j|\!<\! 1, \;j\!=\!1,\cdots,n,\\
{\bm{0}}\; (\in \mC^n)\; \textrm{ otherwise}.
\end{array}}}
$$
Then $\mathbf{\Phi}$ is a continuous map from $\mC^n$ into $\mC^n$. We have that  $\mathbf{\Phi}$ vanishes outside $ \mD^n$, and so
$$
\textstyle 
 \max\limits_{{\bm{z}}\;\! \in\;\! \mD^n }\|\mathbf{\Phi}({\bm{z}})\|_2
= \sup\limits_{{\bm{z}}\;\! \in \;\! \mC^n} \|\mathbf{\Phi}({\bm{z}})\|_2,
$$
where $\|\cdot\|_2$ denotes the usual Euclidean norm in $\mC^n$. 
This implies that there must exist an $r\geq 1$ such that $\mathbf{\Phi}$
maps $K:=r \overline{\mD}^n$ into $K$. As $K$ is compact and convex,
by Brouwer's Fixed Point Theorem (see, e.g., \cite[Theorem~5.28]{Rud}) it follows that there exists a
${\bm{z}}_*\in K$ such that
 $
\mathbf{\Phi}({\bm{z}}_*)={\bm{z}}_*.
$ 
Since $\mathbf{\Phi}$ is zero outside $\mD^n$,  we see that
${\bm{z}}_*\in \mD^n$.  Let ${\bm{z}}_*=(\lambda_1,\cdots,
\lambda_n)$. Then for each $j\in \{1,\cdots, n\}$, we obtain
$$
\textstyle 
\begin{array}{rcl}
\quad \quad \quad 0\!\!\!&=&\!\!\!\lambda_j + \widehat{x}_j(\lambda_1,\cdots, \lambda_n, 
\overline{\lambda_1},\cdots, \overline{\lambda_n},0,\cdots) 
{\scaleobj{0.81}{\prod\limits_{k=1}^n}} \;\!(1\!-\!|\lambda_k|^2)
\\ 
\!\!\!&=&\!\!\! \lambda_j + (\widehat{x}_j \widehat{f}_{n+1})(\lambda_1,\cdots, \lambda_n, 
\overline{\lambda_1},\cdots, \overline{\lambda_n},0,\cdots).\quad \quad \quad \quad \quad (\star\star)
\end{array}
$$
But from ($\star$), we have
$$
\textstyle 
{\scaleobj{0.81}{\sum\limits_{j=1}^n}} (z_j + \widehat{x}_j \widehat{f}_{n+1}) \widehat{y}_j\big|_{(\lambda_1,\cdots, \lambda_n, 
\overline{\lambda_1},\cdots, \overline{\lambda_n},0,\cdots)}= 1 ,
$$
which together with ($\star\star$) yields $0=1$, a contradiction.  As  $n\in \mN$
was arbitrary, it follows that the Bass stable rank of $\partial^{-m}\mathscr{B}$ is
infinite.
\end{proof}

\begin{remarks}$\;$
\begin{enumerate}
\item For Banach algebras, an analogue of the Bass stable rank, called 
the topological stable rank, was introduced  in \cite{Rie}.
  Let $A$ be a commutative complex Banach algebra with unit element
  $1$. The least  $n\in \mN$ for which $U_n(A)$ is dense in $A^n$ is
  called the {\em topological stable rank} of $A$. The {\em topological stable rank of $A$ is infinite} if there is no such $n$. For a commutative unital semisimple complex Banach algebra, the Bass stable
  rank is at  most equal to its topological stable rank (see, e.g., \cite[Corollary~2.4]{Rie}). 
  It follows from Corollary~\ref{Bass_stable_rank_of_A_DS} that the 
  topological stable rank of $\partial^{-m}\mathscr{B}$ is infinite for all $m\in \mN$. 

\item 
  The {\em Krull dimension} of a commutative ring $A$ is the supremum of the
  lengths of chains of distinct proper prime ideals of $A$. If a ring has Krull
  dimension $d$, then its Bass stable rank is at most $d+2$ (see, e.g., 
  \cite{Hei}). It follows from  Corollary~\ref{Bass_stable_rank_of_A_DS} that  the 
  Krull dimension of $\partial^{-m}\mathscr{B}$ is infinite for all $m\in \mN$. 
  \end{enumerate}
\end{remarks}

\end{document}